\documentclass{article}\usepackage{amsmath,amssymb,amsfonts}   
\usepackage{graphicx,bm}

\renewcommand{\v}{\widetilde{v}}
\renewcommand{\u}{\widetilde{u}}
\newcommand{\calU}{{\mathcal U}}
\newcommand{\calG}{{\mathcal G}}
\newcommand{\R}{{\mathbb R}}

\newcommand{\x}{\mathbf{x}}
\newcommand{\y}{\mathbf{y}}
\newcommand{\e}{{\mathrm e}}

\newcommand{\pcb}[1]{\textcolor{blue}{#1}}
\newcommand{\n}{\mathbf n}

\newcommand{\p}{\widetilde{p}}

\newcommand{\hxi}{{\bm \xi}}

\begin{document}

\title{Accumulation time of diffusion in a 3D singularly perturbed domain}

\author{\em 
P. C. Bressloff \\ Department of Mathematics, University of Utah \\
155 South 1400 East, Salt Lake City, UT 84112}

\maketitle

\begin{abstract} 

Boundary value problems for diffusion in singularly perturbed domains (domains with small holes removed from the interior) is a topic of considerable current interest. Applications include intracellular diffusive transport and the spread of pollutants or heat from localized sources. In a previous paper, we introduced a new method for characterizing the approach to steady-state in the case of two-dimensional (2D) diffusion. This was based on a local measure of the relaxation rate known as the accumulation time $T(\x)$. The latter was calculated by solving the diffusion equation in Laplace space using a combination of matched asymptotics and Green's function methods. We thus obtained an asymptotic expansion of $T(\x)$
 in powers of $\nu=-1/\ln \epsilon$, where $\epsilon$ specifies the relative size of the holes. In this paper, we develop the corresponding theory for three-dimensional (3D) diffusion. The analysis is a non-trivial extension of the 2D case due to differences in the singular nature of the Laplace transformed Green's function. In particular, the asymptotic expansion of the solution of the 3D diffusion equation in Laplace space involves terms of order $O((\epsilon/s)^n)$, where $s$ is the Laplace variable. These $s$-singularities have to be removed by partial series resummations in order to obtain an asymptotic expansion of $T(\x)$ in powers of $\epsilon$.

\end{abstract}

\section{Introduction}

There is considerable current interest in solving boundary value problems (BVPs) for two-dimensional (2D) and three-dimensional (3D) diffusion in singularly perturbed domains, where small holes or perforations are removed from the interior \cite{Ward93,Ward93a,Straube07,Bressloff08,Coombs09,Cheviakov11,Chevalier11,Ward15,Coombs15,Bressloff15,Lindsay16,Lindsay17,Grebenkov20,Bressloff21A,Bressloff21B}. Applications range from modeling intracellular diffusion, where interior holes could represent subcellular structures such as organelles or biochemical substrates, to tracking the spread of chemical pollutants or heat from localized sources.

\noindent Roughly speaking, one can divide the various BVPs into two distinct groups. The first treats the holes as totally or partially absorbing traps, and the main focus is determining the first passage time or splitting probability for a single particle to be captured by an interior trap (narrow capture). The second treats the holes as localized sources or reflecting obstacles, and now one is interested in calculating the steady-state solution (if it exists) and the rate of approach to steady state. Both types of BVP can be solved using a combination of matched asymptotic analysis and Green's function methods. This involves obtaining an inner or local solution of the diffusion equation that is valid in a small neighborhood of each hole, and then matching to an outer or global solution that is valid away from each neighborhood. The matching requires taking into account the singular nature of the associated Green's function. However, the details of the matched asymptotic analysis in 2D and 3D domains differ considerably due to corresponding differences in the Green's function singularities. That is, as $|\x-\x_0|\rightarrow 0$, 
\begin{equation}
\label{sing}
G(\x|\x_0)\rightarrow -\frac{1}{2\pi D}\ln|\x-\x_0| \mbox{ in 2D }  \quad G(\x|\x_0)\rightarrow \frac{1}{4\pi D|\x-\x_0|} \mbox{ in 3D}.
\end{equation}
Consequently, an asymptotic expansion of the solution to a BVP in 3D is in powers of $\epsilon$, where $\epsilon$ represents the size of a hole relative to the size of the bulk domain. On the other hand, an analogous expansion in 2D is in powers of $\nu=-1/\ln \epsilon$ at $O(1)$ in $\epsilon$. The slower convergence of $\nu$ in the limit $\epsilon \rightarrow 0$ can be dealt with by summing the logarithmic terms non-perturbatively \cite{Ward93,Ward93a}.

In a recent paper \cite{Bressloff22}, we introduced and analyzed a new quantity for characterizing the rate of relaxation to steady-state in a 2D singularly perturbed domain containing circular holes, based on the so-called accumulation time. The latter is a local measure of the rate of relaxation that has been used extensively within the context of diffusion-based morphogenesis \cite{Berez10,Berez11,Gordon11,Bressloff19}. (Previous studies of singularly perturbed BVPs have considered a global measure of the relaxation rate that is identified with the principal eigenvalue of the Laplacian \cite{Ward93,Ward93a,Coombs09}.) The accumulation time was calculated by solving the diffusion equation in Laplace space, which yielded an asymptotic expansion of the accumulation time
 in powers of $\nu$. In this paper, we develop the corresponding theory for diffusion in 3D singularly perturbed domains containing spherical holes. The analysis is a non-trivial extension of the 2D case due to differences in the singular nature of the Laplace transformed Green's function with respect to the limits $\x\rightarrow \x_0$ and $s\rightarrow 0$, where $s$ is the Laplace variable. In particular, the asymptotic expansion of the solution of the 3D diffusion equation in Laplace space involves terms of order $O((\epsilon/s)^n)$. These $s$-singularities have to be removed by partial series resummations in order to obtain an asymptotic expansion of the accumulation time in powers of $\epsilon$. Surprisingly, in spite of significant differences in the analyses, we find that the $O(1/\epsilon)$ and $O(1) $ contributions to the accumulation time are formally identical to the corresponding terms in 2D under the mappings (from 3D to 2D)
$4\pi D \rightarrow 2\pi D$ and $\epsilon \ell_j \rightarrow \nu_j\equiv -1/{\ln \epsilon \ell_j}$, where $\epsilon \ell_j$ is the radius of the $j$th hole.

The structure of the paper is as follows. In section 2 we  formulate the general problem of diffusion in a 3D singularly-perturbed domain $\Omega$ and define the associated accumulation time  in terms of the Laplace transform of the concentration. The accumulation time is calculated in section 3 by solving the diffusion equation in Laplace space using a combination of matched asymptotic analysis and Green's function methods. Our results are compared with those previously obtained in 2D. In section 4, we relate our analysis to an alternative approach based on an eigenfunction expansion. Finally, in section 5, we illustrate the theory by considering holes in a spherical domain, for which the associated Green's function is known explicitly.

\section{Accumulation time of diffusion in a 3D singularly perturbed domain}

\begin{figure}[b!]
\centering
\includegraphics[width=10cm]{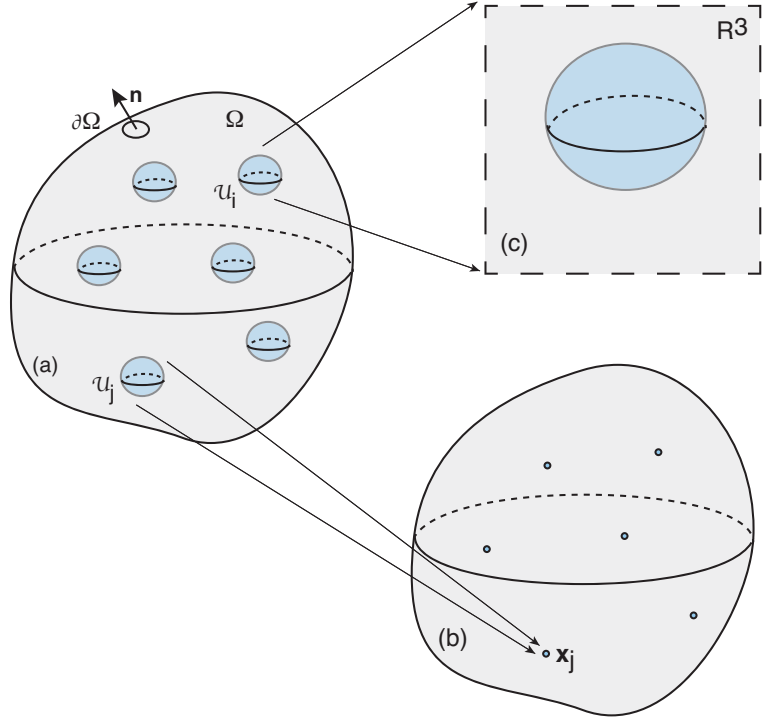} 
\caption{Diffusion in a 3D singularly perturbed domain. (a) Particles diffuse in a bounded domain $\Omega$ containing $N$ small interior holes or perforations denoted by $\calU_j$, $j=1,\ldots,N$. The exterior boundary $\partial \Omega$ is reflecting, whereas $u=\Phi_j$ on the $j$-th interior boundary $\partial \calU_i$. (b) Construction of the outer solution. Each hole is shrunk to a single point. The outer solution can be expressed in terms of the corresponding modified Neumann Green's function and then matched with the inner solution around each hole. (c) Construction of the inner solution in terms of stretched coordinates $\y=\epsilon^{-1}(\x-{\x}_i)$, where ${\x}_i$ is the center of the $i$-th hole. The rescaled radius is $\rho_i=\ell_i$ and the region outside the hole is taken to be $\R^3$ rather than the bounded domain $\Omega$. }
\label{fig1}
\end{figure}

Consider the diffusion equation in a bounded domain $\Omega\subset \R^3$, that is perforated by a set of $N$ small holes denoted by $\calU_k\subset \Omega$, $k=1,\ldots,N$, see Fig. \ref{fig1}(a). The volume of each hole is taken to be $|\calU_j|\sim \epsilon^3 |\Omega|$, $0<\epsilon \ll 1$, with $\calU_j\rightarrow \x_j\in \Omega$ uniformly as $\epsilon \rightarrow 0$, $j=1,\ldots,N$. In addition, the holes are assumed to be well separated with $|\x_i-\x_j|=O(1)$, $j\neq i$, and $\mbox{dist}(x_j,\partial \Omega)=O(1)$ for all $j=1,\ldots,N$. For simplicity, we take each hole to be a sphere with $|\x-\x_j|=\epsilon \ell_j$. We impose a Neumann boundary condition on the external boundary $\partial \Omega$ and inhomogeneous Dirichlet boundary conditions on the interior boundaries $\partial \calU_j$. Let $u(\x,t)$ denote the concentration of freely diffusing particles for $\x\in \Omega\backslash \calU_a$, and $\calU_a\equiv\bigcup_{j=1}^N \calU_j$. Then
\begin{subequations} 
\label{RD1}
\begin{align}
	\frac{\partial u(\x,t)}{\partial t} &= D\nabla^2 u(\x,t),\ \x\in \Omega\backslash \calU_a,	\end{align}
together with the boundary conditions
\begin{equation}
\nabla u(\x,t)\cdot \n=0,\ \x \in \partial \Omega ; \quad u(\x,t)=\Phi_j, \ \x\in \partial \calU_j.
\end{equation}
Here $\n$ is the outward unit normal at a point on $\partial \Omega$. Finally, we impose the initial condition 
\begin{equation}
u(\x,0)=\Gamma_0\delta(\x-\x_0)
\end{equation}
 for some $\x_0 \in  \Omega\backslash \calU_a$, where $\Gamma_0$ is the initial number of molecules introduced into the domain.
 \end{subequations}

Let
\begin{equation}
\label{accu}
Z(\x,t)=1-\frac{u(\x,t)}{u^*(\x)}
\end{equation}
be the fractional deviation of the concentration from steady state. In order to ensure that there is no overshooting (reversal in the sign of $Z(\x,t)$), 
we impose the condition
\begin{equation}
\frac{\Gamma_0}{|\Omega|} < \Phi_j \mbox{ for all }j=1,\ldots,N.
\label{grow}
\end{equation}
Then $1-Z(\x,t)$ represents the fraction of the steady-state concentration that has accumulated at $\x$ by time $t$, and $-\partial_t Z(\x,t)dt$ is the fraction accumulated in the interval $[t,t+dt]$. The accumulation time $T(\x)$ at position $\x$ is then defined as \cite{Berez10,Berez11,Gordon11,Bressloff19}:
\begin{equation}
\label{accu2}
T(\x)=\int_0^{\infty} t\left (-\frac{\partial Z(\x,t)}{\partial t}\right )dt=\int_0^{\infty} Z(\x,t)dt.
\end{equation}
In practice, it is more convenient to calculate the accumulation time in Laplace space. Using the identity 
\[u^*(\x)=\lim_{t\rightarrow \infty} u(\x,t)=\lim_{s\rightarrow 0}s\widetilde{u}(\x,s),\]
where $\u(\x,s)=\int_0^{\infty}\e^{-st}u(\x,t)dt$, and setting $\widetilde{F}(\x,s)=s\widetilde{u}(\x,s)$, the Laplace transform of equation (\ref{accu}) gives
\[s\widetilde{Z}(\x,s)=1-\frac{\widetilde{F}(\x,s)}{\widetilde{F}(\x)},\quad \widetilde{F}(\x)=\lim_{s\rightarrow 0}\widetilde{F}(\x,s)=u^*(\x)\]
and, hence
\begin{eqnarray}
 T(\x)&=&\lim_{s\rightarrow 0} \widetilde{Z}(\x,s) = \lim_{s\rightarrow 0}\frac{1}{s}\left [1-\frac{\widetilde{F}(\x,s)}{\widetilde{F}(\x)}\right ] =-\frac{1}{\widetilde{F}(\x)}
\left .\frac{d}{ds}\widetilde{F}(\x,s)\right |_{s=0}.
\label{acc}
\end{eqnarray}

In light of equation (\ref{acc}), we will work with the diffusion equation (\ref{RD1}) in Laplace space:
\begin{subequations}
\label{masterLT}
\begin{align}
 D\nabla^2\widetilde{u}-s\widetilde{u}  &=-\Gamma_0\delta(\x-\x_0), \quad \x \in \Omega\backslash \calU_a,\\
D\nabla \widetilde{u}(\x,s)\cdot \n&=0,\ \x \in \partial \Omega ,\quad  \widetilde{u}(\x,s)= \frac{\Phi_j}{s},\ \x\in \partial  \calU_j .
\end{align}
\end{subequations}
The Dirac delta function on the right-hand side of equation (\ref{masterLT}a) can be eliminated by introducing the Green's function of the 3D modified Helmholtz equation, 
\begin{align}
\label{GMH}
	D\nabla^2 G(\x,s|\x_0) -sG(\x,s|\x_0) &=-\delta(\x-\x_0) , \ \x\in \Omega,\quad
	 \nabla G(\x,s|x_0)\cdot \n  =0,\ \x \in \partial \Omega, 
	\end{align}
	Two useful features of the Green's function that will play an important role in the subsequent analysis are its singularity structure and its normalization:
	\begin{align}
	\label{propG}
	G(\x,s|\x_0) =\frac{1}{4\pi D|\x-\x_0|}+R(\x,s|\x_0),\quad \int_{\Omega }G(\x,s|\x_0)d\x=\frac{1}{s},
	\end{align}
	where $R(\x,s|\x_0)$ is defined to be the regular part of the Green's function. Finally, taking
\begin{equation}
\label{uv}
\widetilde{u}(\x,s)=\Gamma_0G(\x,s|\x_0)+\v(\x,s),\ \x \in \Omega\backslash \calU_a,
\end{equation}
we have
\begin{subequations} 
\label{masterLT2}
\begin{align}
	D\nabla^2 \v(\x,s ) -s\v(\x,s )&= 0 , \ \x\in \Omega\backslash \calU_a,\\
	 \nabla \v \cdot \n=0, \ \x\in \partial \Omega,\quad  \v&= \frac{\Phi_j}{s}-\Gamma_0G(\x,s|x_0),\ \x\in \partial\calU_i.
	\end{align}
\end{subequations}

\section{Matched asymptotic analysis of the accumulation time in 3D}

The goal of this paper is to derive an asymptotic expansion of the accumulation time (\ref{acc}) in powers of $\epsilon$. We will proceed along analogous lines to studies of the 3D narrow capture problem \cite{Cheviakov11,Coombs15,Bressloff21B}, deriving an inner or local solution of equations (\ref{masterLT}) that is valid in an $O(\epsilon)$ neighborhood of each hole, and then matching to an outer or global solution that is valid away from each neighborhood. However, as previously highlighted in Ref. \cite{Bressloff21B}, the resulting asymptotic expansion of the solution in Laplace space results in terms of order $O((\epsilon/s)^n)$. Therefore, we will have to remove these $s$-singularities in order to obtain an asymptotic expansion of the accumulation time in the limit $s\rightarrow 0$. 

The outer solution is constructed by shrinking each domain $\calU_j$ to a single point $\x_j$, see Fig. \ref{fig1}(b), and expanding according to
\[\u(\x,s)\sim \Gamma_0 G(\x,s|\x_0)+\epsilon\v_1(\x,s)+\epsilon^2 \v_2(\x,s)+\ldots,
\]
where $G$ is the 3D Neumann Green's function, see equation (\ref{GMH}), and
\begin{align}
\label{asym1}
D\nabla^2 \v_n-s\v_n&=0,\, \x\in \Omega\backslash \{\x_1,\ldots,\x_N\};\ \nabla \v_n\cdot \n=0,\, \x\in \partial \Omega .
\end{align}
Equation (\ref{asym1}) is supplemented by a corresponding set of singularity conditions as $\x\rightarrow \x_j$, $j=1,\ldots,N$, which are obtained by matching to the inner solution around each hole. Introducing the stretched local variable ${\mathbf y} = \varepsilon^{-1}(\x-\x_j)$ in a neighborhood of the $j$th hole, see Fig. \ref{fig1}(c), we
set $U(\y,s)=\u(\x_j+\varepsilon \y,s)$ with
\begin{align}
\label{inner}
& D\nabla^2_{\y}U = \epsilon^2 s U,\ |\y| > \ell_j,\quad U(\y,s)=\frac{\Phi_j}{s},\ |\y|=\ell_j.
\end{align}
Substituting the asymptotic expansion 
$U\sim U_0+\epsilon U_1+O(\epsilon^2)$
into (\ref{inner}), we obtain the following pair of equations for the first two terms in the expansion:
\begin{subequations} 
\label{stretch}
\begin{align}
D\nabla_{\y}^2 U_0(\y,s) &=0,\ |\y|>\ell_j;\ U_0(\y,s)= \frac{\Phi_j}{s},\ |\y|=\ell_j,\\ 
	D\nabla_{\y}^2 U_1(\y,s) &=0, \ |\y|>\ell_j;\
	U_1(\y,s)= 0,\ |\y|=\ell_j.
	\end{align}
\end{subequations}
These are supplemented by far-field conditions obtained by matching with the near-field behavior of the outer solution. In order to perform this matching, we need 
to Taylor expand $G(\x,s|\x_0)$ near
the $j$-th target and rewrite it in terms of
stretched coordinates:
\begin{equation}
\label{p0}
G(\x,s|\x_0)  \sim G(\x_j,s|\x_0)+ \epsilon \nabla_{\x} G(\x_j,s|\x_0) \cdot \y +\ldots
\end{equation}

First consider the leading order contribution to the inner solution. Matching the far-field behavior of $U_0$ with the near-field behavior of $\Gamma_0 G(\x,s|\x_0)$ shows that
\begin{equation}
 U_0 \sim \Gamma_0 G(\x_j,s|\x_0) \mbox{ as } |\y|\to \infty. 
\end{equation}
Hence,
\begin{equation}
\label{0U}
U_0 = \frac{\Phi_j}{s} w(\y)+\Gamma_0G(\x_j,s|\x_0)( 1 - w(\y))),
\end{equation}
 with $w(\y)$ satisfying the boundary value problem
\begin{align}
\label{w}
\nabla_{\bf y}^2 w(\y)&=0,\  |\y|>\ell_j ; \quad w(\y)=1,\ |\y|=\ell_j;\ 
w(\y)\rightarrow 0\quad \mbox{as } |\y|\rightarrow \infty. 
\end{align}
In the case of a spherical target of radius $\ell_j$, we have
\begin{equation}
\label{w}
w(\y)=\frac{\ell_j}{|\y|}.
\end{equation}
It now follows that $\v_1$ satisfies equation (\ref{asym1}) together with the singularity condition
\[\v_1(\x,s)\sim \frac{\ell_jV_j(s)}{|\x-\x_j|} \quad \mbox{as } \x\rightarrow \x_j,\]
where
\begin{equation}
\label{Vj}
V_j(s)=\frac{\Phi_j}{s}-\Gamma_0G_{j0}(s), \quad G_{j0}(s)\equiv G(\x_j,s|\x_0).
\end{equation}
Hence,
\begin{equation}
\label{q1}
\v_1(\x,s)= {4\pi}  D\sum_{j=1}^N\ell_jV_j(s)G(\x,s|\x_j).
\end{equation}

The next step is to match the far-field behavior of $U_1$ with the $O(\epsilon)$ term in the expansion of $\Gamma_0 G(\x,s|\x_0)$, see equation (\ref{p0}), together with the non-singular near-field behavior of $\v_1 $ around the $j$-th target. The latter takes the form
\begin{align*}
\v_1(\x,s)&\sim \frac{\ell_jV_j(s)}{|\x-\x_j|}+4\pi D\sum_{i=1}^N\ell_iV_i(s){\mathcal G}_{ij}(s).
\end{align*}
with
\begin{equation}
{\mathcal G}_{ij}(s)=G(\x_i,s|\x_j),\ j\neq i,\quad {\mathcal G}_{jj}(s)=R(\x_j,s|\x_j).
\end{equation}
It follows that
\begin{align}
U_1(\y,s)\rightarrow  \nabla_{\x} G(\x_j,s|\x_0) \cdot \y +4\pi D\sum_{i=1}^N\ell_i V_i(s){\mathcal G}_{ij}(s) \mbox{ as } |\y|\rightarrow \infty .
\end{align}
The first term on the right-hand side generates contributions to the inner solution in the form of first-order spherical harmonics \cite{Bressloff21B}. Since these only affect the outer solution at $O(\epsilon^3)$, we neglect them here.
We thus have
\begin{equation}
\v_1(\x,s)=\chi_j^{(1)}(s) \left (1-\frac{\ell_j}{|\y|}\right )+\mbox{ higher-order harmonics}
\end{equation}
with
\begin{equation}
 \chi_j^{(1)}(s)=4\pi D\sum_{i=1}^N\ell_i V_i(s) {\mathcal G}_{ij}(s) .
 \end{equation}
Finally, $\v_2$ satisfies equation (\ref{asym1}) supplemented by the singularity condition
\[\v_2(\x,s)\sim - \frac{\chi_j^{(1)}(s)\ell_j}{|\x-\x_j|}, \quad \mbox{as } \x\rightarrow \x_j.\]
Using the same steps as in the derivation of $\v_1(\x,s)$, we obtain the result
\begin{equation}
 \v_2(\x,s)= -{4\pi}D\sum_{j=1}^N\ell_j\chi_j^{(1)}(s)G(\x,s|\x_j) .
\end{equation}
In summary, the outer solution has the asymptotic expansion
\begin{equation}
\u(\x,s)\sim \Gamma_0 G(\x,s|\x_0)+4\pi D\sum_{j=1}^N\ell_j\left [\epsilon V_j(s)-\epsilon^2\chi_j^{(1)}(s)+\ldots \right ]G(\x,s|\x_j).
\label{out}
\end{equation}

\subsection{Steady-state solution}

Multiplying equation (\ref{out}) by $s$ and then taking the limit $s\rightarrow 0$ yields the steady-state solution
\begin{align}
\label{split2}
u^*(\x)\sim \Gamma_0 G(\x,s|\x_0)+4\pi D\lim_{s\rightarrow 0}s\sum_{j=1}^N\ell_j\left [\epsilon V_j(s)-\epsilon^2\chi_j^{(1)}(s)+\ldots \right ]G(\x,s|\x_j).
\end{align}
In order to calculate the above limit, we use the result that
\begin{equation}
\label{Gs}
G(\x,s|\x_0)=\frac{1}{s|\Omega|}+G_0(\x,\x_0)+sG_1(\x,\x_0)+O(s^2),
\end{equation}
where $G_0$ is the generalized Neumann Green's function of Laplace's equation:
\begin{subequations}
\label{GM}
\begin{align}
	D\nabla^2 G_0(\x,\x_0)  &=\frac{1}{|\Omega|}-\delta(\x-\x_0) , \ \x\in \Omega,\\
	 \nabla G_0(\x,\x_0)\cdot \n  &=0,\ \x \in \partial \Omega,\quad \int_{\calU}G_0(\x,\x_0)d\x=0, \\
	 G_0(\x,\x_0)&=\frac{1}{4\pi |\x-\x_0|}+R_0(\x,\x_0).
	\end{align}
	\end{subequations}
It follows that the coefficient $V_j$ has the small-$s$ expansion
\begin{equation}
V_j=\frac{\Phi_j-\Gamma_0/|\Omega|}{s}-\Gamma_0G_0(\x_j,\x_0)+O(s).
\end{equation}
Substituting equation (\ref{Gs}) into (\ref{split2}) gives
\begin{align}
\u(\x,s)&\sim\Gamma_0\left [\frac{1}{s|\Omega|}+G_0(\x,\x_0)+O(s)\right ]+4\pi \epsilon D \sum_{j=1}^N\ell_j \bigg \{\left [\frac{\widehat{\Phi}_j}{s}-\Gamma_0\calG_{j0}^{(0)}-s\Gamma_0\calG_{j0}^{(1)}+O(s^2)\right ] \nonumber \\
&\quad -4\pi \epsilon D \sum_{k=1}^N \ell_k \left [\frac{\widehat{\Phi}_k}{s}-\Gamma_0\calG_{k0}^{(0)}-s\Gamma_0\calG_{k0}^{(1)}+O(s^2) \right ]\left [\frac{1}{s|\Omega|}+{\mathcal G}_{kj}^{(0)}+s{\mathcal G}_{kj}^{(1)}+O(s^2)\right ]\bigg \}\nonumber \\
&\quad \times \left [\frac{1}{s|\Omega|}+G_0(\x,\x_j)+sG_1(\x,\x_j)+O(s^2)\right ]+O(\epsilon^3).
\end{align}
We have set $\calG_{k0}^{(n)}=G_n(\x_k,\x_0)$ and $\widehat{\Phi}_j=\Phi_j-{\Gamma_0}/{|\Omega|}$.
Rearranging the various terms and multiplying by $s$ yields the asymptotic expansion
\begin{align}
&s\u(\x,s)\sim\frac{\Gamma_0}{|\Omega|}\bigg \{1-4\pi \epsilon D \sum_{j=1}^N\ell_j \calG_{j0}^{(0)}\bigg \}+4\pi \epsilon D \sum_{j=1}^N\ell_j\widehat{\Phi}_jG_0(\x,\x_j)+O(\epsilon^2)\nonumber \\
&\quad +\frac{4\pi \epsilon D}{s|\Omega|} \bigg \{ \sum_{j=1}^N\ell_j \widehat{\Phi}_j-4\pi \epsilon D \sum_{j,k=1}^N  \ell_j\ell_k\widehat{\Phi}_k\bigg [G_0(\x,\x_j)+\calG_{kj}^{(0)}\bigg ]+\frac{4\pi \epsilon D\Gamma_0}{|\Omega|}\sum_{j=1}^N\ell_j \ell_k\calG_{j0}^{(0)} +O(\epsilon^2)\bigg \}\nonumber \\
&\quad -\left (\frac{4\pi \epsilon D}{s|\Omega|}\right  )^2\bigg \{\sum_{j,k=1}^N\ell_j \ell_k\widehat{\Phi}_k+O(\epsilon)\bigg \}+s\bigg \{\Gamma_0G_0(\x,\x_0) +4\pi \epsilon D \sum_{j=1}^N\ell_j\widehat{\Phi}_jG_1(\x,\x_j)\nonumber \\
&\quad 
-4\pi \epsilon D \Gamma_0 \sum_{j=1}^N \ell_j[\calG_{j0}^{(0)}G_0(\x,\x_j)+\calG_{j0}^{(1)}/|\Omega|]+O(\epsilon^2)\bigg \}. 
\label{JLT2}
\end{align}

The $\epsilon$-expansion in equation (\ref{JLT2}) indicates a potential problem in taking the limit $s\rightarrow 0$. More specifically, there exist terms involving factors of $\epsilon/s$ that will become arbitrarily large in the small-$s$ limit and thus lead to a breakdown of the $\epsilon$ expansion. 
This issue was previously encountered in an analysis of first passage time problems in 3D singularly perturbed domains with small traps \cite{Bressloff21B}. In the latter study, we calculated the Laplace-transformed flux into each trap, which acted as the generator for the first passage time moments in the limit $s\rightarrow 0$. In contrast to the analysis of steady-state problems for diffusion, where the goal is to calculate the outer solution in the bulk of the domain, the focus of narrow capture problems is the inner solution around each trap. Nevertheless, the methods developed in Ref. \cite{Bressloff21B} can be adapted to eliminate the singularities in equation (\ref{JLT2}). That is, we proceed by treating equation (\ref{JLT2}), including higher-order terms, as a triple expansion in $\epsilon$, $s$ and $\Lambda$, with
\begin{equation}
\Lambda =\frac{ 4\pi   \epsilon D \bar{\ell}}{s|\calU|},\quad \bar{\ell} =\sum_{j=1}^N\ell_j.
\end{equation}
This then converts a subset of terms at $O(\epsilon^{n})$ to $O(\epsilon^{r} \Lambda^{n-r})$ terms, $0\leq r\leq n$. At each order of $\epsilon$, we obtain infinite power series in $\Lambda$ that can be summed to remove all singularities in the limit $s\rightarrow 0$. In order to illustrate the basic idea, consider the sum of the first terms on the second and third lines of equation (\ref{JLT2}). Inclusion of higher-order contributions leads to a geometric series in $\Lambda$ that can be summed explicitly:  
\begin{equation}
\label{I1}
{\mathcal I}_1(\Lambda)\equiv \left (\frac{\Lambda}{\bar{\ell}}\sum_{n\geq 0}(-1)^n\Lambda^n\right )\sum_{j=1}^N\ell_j \widehat{\Phi}_j=\frac{1}{\bar{\ell}}\frac{\Lambda}{1+\Lambda}\sum_{j=1}^N\ell_j \widehat{\Phi}_j=\frac{\Lambda}{1+\Lambda}[\overline{\Phi}-\Gamma/|\Omega|].
\end{equation}
Similarly, combining the last term on the second line of (\ref{JLT2}) with higher-order contributions
\begin{equation}
\label{I2}
\epsilon {\mathcal I}_2(\Lambda)\equiv  \left (\Lambda \sum_{n\geq 0}(-1)^n\Lambda^n\right )\frac{4\pi \epsilon D\Gamma_0}{|\Omega|} \sum_{j=1}^N\ell_j \calG_{j0}^{(0)}=\frac{\Lambda}{1+\Lambda}\frac{4\pi \epsilon D\Gamma_0}{|\Omega|} \sum_{j=1}^N\ell_j \calG_{j0}^{(0)}.
\end{equation}
Finally, combining the middle term of the second line of equation (\ref{JLT2}) with high-order contributions yields
\begin{align}
&\epsilon {\mathcal I}_3(\Lambda)\equiv -\left (\Lambda \sum_{n\geq 0}(-1)^n\Lambda^n\right )\frac{4\pi \epsilon D}{\bar{\ell}} \sum_{j,k=1}^N  \ell_j\ell_k\widehat{\Phi}_k\bigg [G_0(\x,\x_j)+\calG_{kj}^{(0)}\bigg ]\nonumber \\
&\quad +\Lambda^2\sum_{m\geq 0}(m+1)(-\Lambda)^m\frac{4\pi \epsilon D}{\bar{\ell}^2} \sum_{i,j,k=1}^N \ell_i \ell_j\ell_k\widehat{\Phi}_i\calG_{kj}^{(0)}
\label{I3}\\
&=-\frac{\Lambda}{1+\Lambda}\frac{4\pi \epsilon D}{\bar{\ell}} \sum_{j,k=1}^N  \ell_j\ell_k\widehat{\Phi}_k\bigg [G_0(\x,\x_j)+\calG_{kj}^{(0)}\bigg ]+\frac{\Lambda^2}{(1+\Lambda)^2}
\frac{4\pi \epsilon D}{\bar{\ell}} [\overline{\Phi}-\Gamma/|\Omega|]\sum_{j,k=1}^N  \ell_j\ell_k\calG_{kj}^{(0)}.\nonumber 
\end{align}
We have used the following results for geometric series:
\begin{subequations}
\label{insum}
\begin{align}
 \Lambda\sum_{m\geq 0}(-\Lambda)^m&=\Lambda(1-\Lambda+\Lambda^2 \ldots )=\frac{\Lambda}{1+\Lambda},\\
\Lambda^2\sum_{m\geq 0}(m+1)(-\Lambda)^m&=\Lambda^2(1-2\Lambda+3\Lambda^2\ldots )=\Lambda^2 \frac{d}{d\Lambda}\frac{\Lambda}{1+\Lambda}=\frac{\Lambda^2}{(1+\Lambda)^2} .
\end{align}
\end{subequations}

Having performed the various partial summations, equation (\ref{JLT2}) can be rewritten in the more compact form
\begin{align}
s\u(\x,s)&\sim\frac{\Gamma_0}{|\Omega|}\bigg \{1-4\pi \epsilon D \sum_{j=1}^N\ell_j \calG_{j0}^{(0)}\bigg \}+4\pi \epsilon D \sum_{j=1}^N\ell_j\widehat{\Phi}_jG_0(\x,\x_j)\nonumber \\&\quad +{\mathcal I}_1(\Lambda)+\epsilon{\mathcal I}_2(\Lambda)+\epsilon {\mathcal I}_3(\Lambda)+s\bigg \{\Gamma_0G_0(\x,\x_0) +4\pi \epsilon D \sum_{j=1}^N\ell_j\widehat{\Phi}_jG_1(\x,\x_j)\nonumber \\
&\quad 
-4\pi \epsilon D \Gamma_0 \sum_{j=1}^N \ell_j[\calG_{j0}^{(0)}G_0(\x,\x_j)+\calG_{j0}^{(1)}/|\Omega|]\bigg \}+O(\epsilon^2,s^2). 
\label{JLT2a}
\end{align}
We can now safely take the limit $s\rightarrow 0$ for fixed $\epsilon>0$ in equations (\ref{I1})--(\ref{I3}), since $\Lambda \rightarrow \infty$ and $\Lambda/(1+\Lambda)\rightarrow 1$. We thus obtain the following asymptotic expansion of the steady state to $O(\epsilon)$:
\begin{align}
u^*(\x)&=\lim_{s\rightarrow 0}s\u(\x,s)\sim\frac{\Gamma_0}{|\Omega|}+\frac{1}{\bar{\ell}}\sum_{j=1}^N\ell_j \widehat{\Phi}_j+4\pi \epsilon D \sum_{j=1}^N\ell_j\widehat{\Phi}_jG_0(\x,\x_j)\\
&\quad -\frac{1}{\bar{\ell}} \bigg \{4\pi \epsilon D \sum_{j,k=1}^N  \ell_j\ell_k\widehat{\Phi}_k\bigg [G_0(\x,\x_j)+\calG_{kj}^{(0)}\bigg ]\bigg \}+\frac{4\pi \epsilon D}{\bar{\ell}} \frac{1}{\bar{\ell}}\sum_{i=1}^N \ell_i\widehat{\Phi}_i\sum_{j,k=1}^N  \ell_j\ell_k\calG_{kj}^{(0)} .\nonumber 
\end{align}
Finally, noting that $\widehat{\Phi}_j=\Phi_j-\Gamma_0/|\Omega|$, we see that any dependence on the initial distribution $u(\x,0)=\Gamma_0\delta(\x-x_0)$ vanishes and 
\begin{align}
u^*(\x)&\sim\overline{\Phi}+4\pi \epsilon D \sum_{j=1}^N\ell_j[\Phi_j-\overline{\Phi}]G_0(\x,\x_j)-\frac{4\pi \epsilon D}{\bar{\ell}} \sum_{j,k=1}^N  \ell_j\ell_k[{\Phi}_k-\overline{\Phi}] \calG_{kj}^{(0)}+O(\epsilon^2)
\label{JLT3}
\end{align}
with $\overline{\Phi}={\bar{\ell}}^{-1}\sum_{j=1}^N\ell_j  {\Phi}_j$.
Equation (\ref{JLT3}) is identical to the result obtained by directly solving the steady-state diffusion equation using matched asymptotics \cite{Cheviakov11,Coombs15}. However, the advantage of working in Laplace space is that one can also calculate the accumulation time.

\subsection{Accumulation time}

In order to calculate the accumulation time according to equation (\ref{acc}), we need to differentiate both sides of equation (\ref{JLT2a}) with respect to $s$. 
\begin{align}
\frac{ds\u(\x,s)}{ds}&\sim\bigg \{\Gamma_0G_0(\x,s|\x_0) +4\pi \epsilon D \sum_{j=1}^N\ell_j\widehat{\Phi}_jG_1(\x,\x_j)-4\pi \epsilon D \Gamma_0 \sum_{j=1}^N \ell_j\calG_{j0}^{(0)}G_0(\x,\x_j)\nonumber \\
&\quad 
+\calG_{j0}^{(1)}/|\Omega|]\bigg \} +\frac{d\Lambda}{ds} \frac{d}{d\Lambda} \left ({\mathcal I}_1(\Lambda)+\epsilon{\mathcal I}_2(\Lambda)+\epsilon {\mathcal I}_3(\Lambda)\right )+O(\epsilon^2,s).
\label{JLT3a}
\end{align}
Using the results
\[\frac{d\Lambda}{ds}=-\frac{1}{s^2} \frac{4\pi \epsilon D\bar{\ell}}{|\Omega|}=-\frac{\Lambda}{s}, \quad \frac{d}{d\Lambda}\frac{\Lambda^n}{(1+\Lambda)^n}=-\frac{n{\Lambda}^{n-1}}{(1+{\Lambda})^{n+1}},
\]
it follows that
\begin{equation}
\frac{d}{ds}\frac{\Lambda^n}{(1+\Lambda)^n}=-\frac{1}{s}\frac{n{\Lambda}^n}{(1+{\Lambda})^{n+1}}\rightarrow -\frac{n|\Omega|}{4\pi \epsilon D\bar{\ell}} \mbox{ as } \ {s\rightarrow 0}.
\end{equation}
Therefore, using equations (\ref{I1})--(\ref{I3}) we have
\begin{align}
{\mathcal F}(\x)&\equiv \left .\frac{d\widetilde{F}(\x,s)}{ds}\right |_{s=0}=\frac{1}{4\pi \epsilon D\bar{\ell}}[\Gamma_0-|\Omega|\overline{\Phi}]+\Gamma_0G_0(\x,\x_0)-\frac{\Gamma_0}{\bar{\ell}} \sum_{j=1}^N\ell_j \calG_{j0}^{(0)}\nonumber \\
&\quad +\frac{|\Omega|}{\bar{\ell}^2} \sum_{j,k=1}^N  \ell_j\ell_k\widehat{\Phi}_k\bigg [G_0(\x,\x_j)+\calG_{kj}^{(0)}\bigg ]-
\frac{2}{\bar{\ell}^2} [|\Omega|\overline{\Phi}-\Gamma_0]\sum_{j,k=1}^N  \ell_j\ell_k\calG_{kj}^{(0)}+O(\epsilon)\nonumber \\
&=\frac{1}{4\pi \epsilon D\bar{\ell}}[\Gamma_0-|\Omega|\overline{\Phi}]+\Gamma_0G_0(\x,\x_0)-\frac{\Gamma_0}{\bar{\ell}} \sum_{j=1}^N\ell_j G_0(\x_j,\x_0)-\frac{\Gamma_0-|\Omega|\overline{\Phi}}{\bar{\ell}}\sum_{j=1}^N\ell_jG_0(\x,\x_j)\\
&\quad +\frac{\Gamma_0-|\Omega|\overline{\Phi}}{\bar{\ell}^2} \sum_{j,k=1}^N  \ell_j\ell_k\calG_{kj}^{(0)} +\frac{|\Omega|}{\bar{\ell}^2} \sum_{j,k=1}^N  \ell_j\ell_k[\widehat{\Phi}_k-\overline{\Phi}]\calG_{kj}^{(0)} +O(\epsilon).\nonumber
\end{align}
Since $\widehat{\Phi}_k=\Phi_k-\Gamma_0/|\Omega|$, it follows that
\begin{align}
{\mathcal F}(\x)\sim \frac{\Gamma_0-|\Omega|\overline{\Phi}}{4\pi \epsilon D\bar{\ell}}+{\mathcal F}_0(\x)+O(\epsilon)
\end{align}
with
\begin{align}
\label{F03D}
{\mathcal F}_0(\x)
&=\Gamma_0G_0(\x,\x_0)-\frac{\Gamma_0}{\bar{\ell}} \sum_{j=1}^N\ell_j G_0(\x_j,\x_0)-\frac{\Gamma_0-|\Omega|\overline{\Phi}}{\bar{\ell}}\sum_{j=1}^N\ell_jG_0(\x,\x_j)\\
&\quad +\frac{\Gamma_0-|\Omega|\overline{\Phi}}{\bar{\ell}^2} \sum_{j,k=1}^N  \ell_j\ell_k\calG_{kj}^{(0)} +\frac{|\Omega|}{\bar{\ell}^2} \sum_{j,k=1}^N  \ell_j\ell_k[\widehat{\Phi}_k-\overline{\Phi}]\calG_{kj}^{(0)} +O(\epsilon).\nonumber
\end{align}

Finally, substituting for ${\mathcal F}(\x)$ into equation (\ref{acc}) and using equation (\ref{JLT3}) for $u^*(\x)$ yields the following result for the accumulation time for diffusion in a 3D singularly perturbed domain:
\begin{align}
T(\x)&=\frac{|\Omega|\overline{\Phi}-\Gamma_0}{4\pi \epsilon D\bar{\ell}\overline{\Phi}}-\frac{{\mathcal F}_0(\x)}{\overline{\Phi}}+\frac{|\Omega|\overline{\Phi}-\Gamma_0}{  \bar{\ell}\overline{\Phi}^2}\left [\sum_{j=1}^N\ell_j[\Phi_j-\overline{\Phi}]G_0(\x,\x_j)-\frac{1}{\bar{\ell}} \sum_{j,k=1}^N  \ell_j\ell_k[{\Phi}_k-\overline{\Phi}] \calG_{kj}^{(0)}\right ]\nonumber \\
&\qquad +O(\epsilon).
\label{Tres3D}
\end{align}
The expression for the accumulation time simplifies considerably in the case of $N$ identical interior boundary conditions, $\Phi_j= {\Phi}$ and identical hole sizes $\ell_j={\ell}$, $j=1,\ldots,N$:
\begin{align}
\label{one}
T(\x)&= \frac{|\Omega|   {\Phi}-\Gamma_0}{4\pi \epsilon D N \ell {\Phi}}-\frac{\Gamma_0}{ {\Phi}}\left [ G_0(\x,\x_0)-\frac{1}{N}\sum_{j=1}^NG_0(\x_j,\x_0)\right ]\\
&\quad -\frac{|\Omega|{\Phi}-\Gamma_0}{ N  {\Phi}}\bigg \{  \sum_{j=1}^NG_0(\x,\x_j)-\frac{1}{N}\sum_{i,j=1}^N\calG_{ij}^{(0)}\bigg \}+O(\epsilon).\nonumber 
\end{align}
In either case, the leading order contribution to $T(\x)$ is the constant 
\begin{equation}
\mu_0\equiv \frac{|\Omega|\overline{\Phi}-\Gamma_0}{4\pi  \epsilon N\overline{\Phi}}.
\end{equation}
Note that $\mu_0>0$ due to the condition (\ref{grow}), which ensures that the accumulation time is positive.
Moreover, $T(\x)\rightarrow \infty$ as $\epsilon \rightarrow 0$. This singular behavior as the size of the holes shrinks to zero is related to the fact that
$\lim_{\epsilon \rightarrow 0}u^*(\x)=\overline{\Phi}$,
whereas the steady-state in the absence of any holes is $\Gamma_0/|\Omega|$. In other words, the limits $\epsilon\rightarrow 0$ and $t\rightarrow \infty$ do not commute.

\subsection{Comparison with the accumulation time for 2D diffusion}

In our previous paper \cite{Bressloff22},  we developed an analogous asymptotic analysis of the accumulation time $T(\x)$ for diffusion in 2D singularly perturbed domains. However, the details of the matched asymptotic analysis differed considerably from the 3D case, reflecting differences in the singular nature of the modified Helmholtz Green's function, see equation (\ref{sing}). Consequently, in 2D we obtained an asymptotic expansion of $T(\x)$ in powers of $\nu=-1/\ln \epsilon$ at $O(1)$ in $\epsilon$. On the other hand, taking the small-$s$ limit was relatively straightforward. 
Surprisingly, in spite of significant differences in the analyses, the $O(1/\epsilon)$ and $O(1) $ contribution to $T(\x)$ in equation (\ref{Tres3D}) are formally identical to the corresponding terms in 2D under the mappings (from 3D to 2D)
\[4\pi D \rightarrow 2\pi D,\quad \epsilon \ell_j \rightarrow \nu_j\equiv -\frac{1}{\ln \epsilon \ell_j},\quad \epsilon \bar{\ell}\rightarrow \sum_{j=1}^N\nu_j ,\]
see equations (4.34)--(4.37) of Ref. \cite{Bressloff22}.

As originally shown by Ward and Keller \cite{Ward93} within the context of 2D and 3D eigenvalue problems, it is possible to generalize the asymptotic analysis of the accumulation time to more general hole shapes such as ellipsoids by applying classical results from electrostatics. For example, given a general shape $\calU_j \subset \R^3$, the solution to equation (\ref{stretch}a) is given by equation (\ref{0U}) with $w(\y)$ having the far-field behavior
\begin{equation}
w(\y)\sim \frac{C_j}{|\y|}+\frac{{\bf P}_j\cdot \y}{|\y|^3}+\ldots \mbox{as } |\y|\rightarrow \infty.
\end{equation}
Here $C_j$ is the capacitance and ${\bf P}_j$ the dipole vector of an equivalent charged conductor with the shape $\calU_j$. (For a sphere, $C_j=\ell_j$ and ${\bf P}_j=0$). It turns out that the $O(\epsilon)$ and $O(\epsilon^2)$ contributions to the accumulation time only depend on $C_j$ so that equation (\ref{Tres3D}) still holds on making the replacements $\ell_j\rightarrow C_j$ for $j=1,\ldots,N$. Similarly, in 2D one simply sets $\nu_j=-1/\ln \epsilon d_j$ with $d_j$ the associated logarithmic capacitance.

\section{Eigenfunction expansion}

Characterizing the relaxation to steady state in terms of the $\x$-dependent accumulation time $T(\x)$ is significantly different from the standard method based on an eigenvalue expansion \cite{Ward93,Ward93a,Coombs09}. Consider the set of eigenpairs of the negative Laplacian in the given singularly perturbed domain, which are denoted by $(\lambda_n,\phi_n(\x))$ for $n\geq 0$ with $0<\lambda_0<\lambda_1 <\lambda_2\ldots$ and 
\[\int_{\Omega\backslash \calU_a}\phi_n(\x)\phi_m(\x)d\x=\delta_{n,m}.\]
 Then
\begin{equation}
\label{eig2}
u(\x,t)-u^*(\x)=\sum_{n\geq 0}c_n \phi_n(\x)\e^{-\lambda_n t}\approx c_0\phi_0(\x)\e^{-\lambda_0 t},
\end{equation}
where $\lambda_0$ is the smallest nonzero eigenvalue. Since $\lambda_0=O(\epsilon)$, it can be calculated by solving the singularly perturbed BVP \cite{Ward93,Cheviakov11} 
 \begin{subequations}
 \label{BVPeig}
\begin{align}
&D\nabla^2 \phi_0+\lambda_0\phi_0=0,\ \x\in \Omega\backslash \calU_a,\quad \nabla\phi_0 \cdot \n =0,\ \x \in \partial \Omega, \ \int_{\Omega\backslash \calU_a}\phi_0^2(\x)d\x=1,\\
&\phi_0=0,\ \x\in \partial \calU_j,\quad j=1,\ldots,N,
\end{align}
\end{subequations}
Following \cite{Cheviakov11}, we expand the principal eigenvalue as
\begin{equation}
\lambda_0 =\epsilon \lambda_0^{(1)}+\epsilon^2\lambda_0^{(2)}+\ldots
\end{equation}
Similarly, the outer eigenfunction is expanded as  
\begin{equation}
\phi_0=\phi_0^{(0)}+\epsilon \phi_0^{(1)}+\epsilon^2 \phi_0^{(2)}+\ldots,
\end{equation}
where $\phi_0^{(0)}=|\Omega|^{-1/2}$. In particular,
\begin{subequations}
\label{eigen}
\begin{align}
&D\nabla^2 \phi_0^{(1)}=-\lambda_0^{(1)}\phi_0^{(0)},\ \x\in \Omega\backslash \{\x_1,\ldots,\x_N\},\quad \nabla \phi_0^{(1)}\cdot \n =0,\ \x \in \partial \Omega, \ \int_{\Omega}\phi_0^{(1)}(\x)d\x=0\\
&D\nabla^2 \phi_0^{(2)}=-\lambda_0^{(2)}\phi_0^{(0)}-\lambda_0^{(1)} \phi_0^{(1)},\ \x\in \Omega\backslash \{\x_1,\ldots,\x_N\},\quad \nabla \phi_0^{(2)}\cdot \n =0,\ \x \in \partial \Omega, \nonumber \\
&\quad  \int_{\Omega}\phi_0^{(2)}(\x)d\x=-\frac{1}{2\phi_0^{(0)}}\int_{\Omega}[\phi_0^{(1)}(\x)]^2d\x.
\end{align}
\end{subequations}
The matching of $\phi_0^{(1)}$ and $\phi_0^{(2)}$ with the inner solution around each hole will yield singularity conditions as $\x\rightarrow \x_j$, $j=1,\ldots,N$. The inner eigensolution is expanded as
$U=U_0+\epsilon U_1+\ldots$
with
\begin{align}
\nabla_{\y}^2 U_k =0, \ |\y|>\ell_j,\quad U_k(\y)=0,\ |\y|=\ell_j
\end{align}
for $k\leq 2$. The near-field behavior of the outer eigenfunction as $\x\rightarrow \x_j$ has to match the far-field behavior of the inner solution as $\y=\epsilon^{-1}|\x-\x_j|\rightarrow \infty$.

The first matching condition is $U_0\rightarrow \phi_0^{(0)}$ as $|\y|\rightarrow \infty$, which means that
\begin{equation}
U_0=\phi_0^{(0)}(1-w(\y)) ,\quad w(\y)=\frac{\ell_j}{|\y|}.
\end{equation}
The singularity condition for $\phi_0^{(1)}$ is thus $\phi_0^{(1)}\sim -\phi_0^{(0)}\ell_j/|\x-\x_j|$ as $\x\rightarrow \x_j$, which implies that
\begin{equation}
\phi_0^{(1)}=-4\pi D \phi_0^{(0)} \sum_{j=1}^N \ell_j G_0(\x,\x_j),
\end{equation}
where $G_0$ is the 3D Green's function satisfying equation (\ref{GM}). Requiring that the solution for $\phi_0^{(1)}$ satisfies equation (\ref{eigen}a) yields
\begin{equation}
\lambda_0^{(1)}=\frac{4\pi D}{|\Omega|} \sum_{j=1}^N\ell_j.
\end{equation}
The near-field behavior of $\phi_0^{(1)}$ is
\begin{equation}
\phi_0^{(1)}\sim -\frac{\ell_j \phi_0^{(0)}}{|\x-\x_j|} -4\pi D \phi_0^{(0)} \sum_{k=1}^N \ell_k\calG_{jk}^{(0)}.
\end{equation}
Matching with the far-field behavior of $U_1$ gives
\begin{equation}
U_1=-4\pi D \phi_0^{(0)} \sum_{k=1}^N \ell_k\calG_{jk}^{(0)}\left (1-\frac{\ell_j}{|\y|}\right ).
\end{equation}
It follows that the singular behavior of $\phi_0^{(2)}$ is
\begin{equation}
\phi_0^{(2)} \sim - \phi_0^{(0)} \frac{ \ell_j\chi_j}{|\x-\x_j|}\mbox{ as } \x\rightarrow \x_j,\quad \chi_j=4\pi D \sum_{k=1}^N  \ell_k\calG_{jk}^{(0)},
\end{equation}
which means that $\phi_0^{(2)}$ satisfies the equation
\begin{align}
&D\nabla^2 \phi_0^{(2)}=-\lambda_0^{(2)}\phi_0^{(0)}-\lambda_0^{(1)} \phi_0^{(1)}-4\pi D\sum_{j=1}^N \ell_j \chi_j \delta(\x-\x_j),\ \x\in \Omega,\quad \nabla \phi_0^{(2)}\cdot \n =0,\ \x \in \partial \Omega.
\end{align}
Applying the divergence theorem with $\int_{\Omega}\phi_0^{(1)}(\x)d\x=0$ yields
\begin{equation}
\lambda_0^{(2)}=-\frac{4\pi D}{|\Omega|} \sum_{j=1}^N\ell_j\chi_j .
\end{equation}

Hence, the principal eigenvalue is given by
\begin{equation}
\lambda_0 = \frac{4\pi \epsilon D}{|\Omega|} \left (\sum_{j=1}^N\ell_j- {4\pi \epsilon D }\sum_{j,k=1}^N \ell_j\ell_k\calG_{jk}^{(0)} \right )+O(\epsilon^3).
\end{equation}
The inverse of the principal eigenvalue can be identified as a global measure of the relaxation rate \cite{Ward93,Cheviakov11}:
\begin{equation}
\tau_0\equiv\frac{1}{\lambda_0} = \frac{|\Omega|}{4\pi \epsilon D \bar{\ell}} \left (1+ \frac{4\pi \epsilon D }{\bar{\ell}}\sum_{j,k=1}^N \ell_j\ell_k\calG_{jk}^{(0)} \right )+O(\epsilon).
\end{equation}
Note that, in contrast to the accumulation time, $\tau_0$ is independent of the initial density and the boundary values $\Phi_j$. As expected, $\tau_0\rightarrow \infty$ as $\epsilon \rightarrow 0$. 
One important advantage of the accumulation time, beyond the fact that it includes local information about the relaxation process, is that it can be calculated without recourse to a spectral decomposition, and thus does nor rely on the existence of a sufficiently large spectral gap. However, one could use the eigenfunction expansion to obtain an approximation of the accumulation time in terms of the principal eigenvalue and eigenfunction. That is, substituting equation (\ref{eig2}) into (\ref{accu2}) implies that
\begin{equation}
T(x)=-\frac{1}{u^*(x)}\int_0^{\infty} \sum_{n\geq 0}\phi_n(x)\e^{-\lambda_n t}dt=-\sum_{n=0}^{\infty}\frac{c_n\phi_n(x)}{\lambda_n u^*(x)},
\end{equation}
which is non-singular since $\lambda_n >0$ for all $n\geq 0$. Keeping only the first term in the series expansion then yields the truncated accumulation time
\begin{equation}
T_0(\x)=-\frac{c_0\phi_0(x)}{\lambda_0 u^*(x)}.
\end{equation}
For simplicity, consider the homogenous case $\Phi_j=\Phi$ and $\ell_j=\ell$ for all $j=1,\ldots,N$, such that $u^*(\x)=\Phi$. The principal eigenvalue and eigenfunction have the asymptotic expansions
\begin{equation}
\lambda_0 = \frac{4\pi \epsilon D N \ell  }{|\Omega|} \left (1- \frac{4\pi \epsilon D \ell}{N}\sum_{j,k=1}^N \calG_{jk}^{(0)} \right )+O(\epsilon^3),
\end{equation}
and
\begin{align}
\phi_0(\x)\sim \frac{1}{\sqrt{|\Omega|}} -\frac{4\pi \epsilon D  \ell  }{\sqrt{|\Omega|}} \sum_{j=1}^N G_0(\x,\x_j)+O(\epsilon^2)
\end{align}

It remains to calculate the coefficient $c_0$. Setting $t=0$ in equation  (\ref{eig2}) gives
\begin{equation}
\Gamma_0\delta(\x-\x_0) -\Phi=\sum_{n=0}^{\infty}c_n\phi_n(\x).
\end{equation}
Multiplying both sides by $\phi_0(\x)$, integrating with respect to $\x$ and imposing orthonormality of the eigenfunctions yields
\begin{equation}
c_0=\Gamma_0\phi_0(\x_0)- \Phi \int_{\Omega} \phi_0(\x)d\x  .
\end{equation}
Substituting the solution for $\phi_0(\x)$ and using the normalization
condition $\int_{\Omega}G_0(\x|\x_0)d\x=0$, we have
\begin{align}
c_0&=\frac{\Gamma_0-|\Omega| \Phi}{\sqrt{|\Omega|}}-\frac{4\pi   \epsilon D  \ell  }{\sqrt{|\Omega|}} \sum_{j=1}^N G_0(\x_0,\x_j)+O(\epsilon^2).
\end{align}
Combining our various results and comparing with equation (\ref{one}) for the full accumulation time shows that
\begin{equation}
T(\x)\sim T_0(\x)-\frac{\Gamma_0}{\Phi}G_0(\x,\x_0)+O(\epsilon)
\end{equation}
It can be seen that the difference between the two is maximized in a neighborhood of the initial position $\x_0$.
Again, we obtained an analogous result for the accumulation time in 2D up to $O(\nu)$ \cite{Bressloff22}.

\section{Examples} 

\begin{figure}[b!]
\centering
\includegraphics[width=12cm]{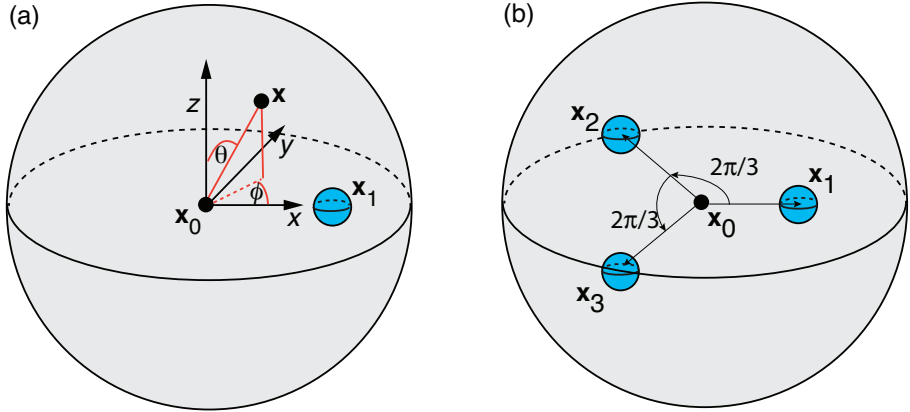} 
\caption{Holes in the unit sphere. (a) Single spherical hole whose center $\x_1$ is located along the $x$-axis of the unit sphere, and the initial concentration is localized at the origin, $\x_0=(0,0,0)$ (b) Triplet of identical spherical holes evenly distributed in the horizontal mid-plane of the unit sphere ($\theta =\pi/2$).}
\label{fig2}
\end{figure}

\begin{figure}[b!]
\centering
\includegraphics[width=13cm]{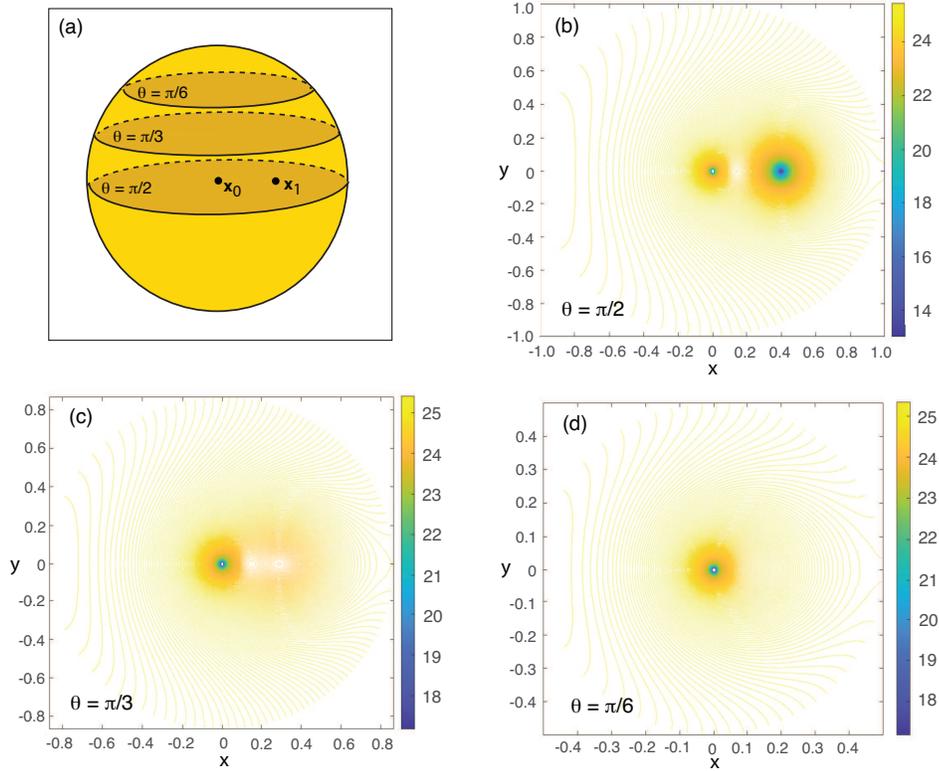} 
\caption{Accumulation time $T(\x)$ in the unit sphere with a single-hole, see Fig. \ref{fig2}(a). The accumulation time is sampled across several horizontal sections of the sphere as indicated in (a). This generates contour plots of $T(\x)$ in the $x-y$ plane for $\x=r(\sin \theta \cos \phi,\sin \theta\sin \phi.\cos \theta) $ with $0\leq r \leq 1$, $0\leq \phi <2\pi$ and fixed $\theta$: (b) $\theta=\pi/2$; (c) $\theta=\pi/3$; (d) $\theta=\pi/6$. Other parameter values are $a=0.4$, $\Gamma_0=1$, $\epsilon=0.01$ and $D=1$.}
\label{fig3}
\end{figure}

\subsection{Single target in the unit sphere}

As our first example, consider the 3D configuration shown in Fig. \ref{fig2}(a). The domain $\Omega$ is taken to be the unit sphere with a single hole placed at $\x_{1}=(a,0,0)$. The boundary condition is
\begin{equation}
\nabla u(\x,t)\cdot {\bf n}=0,\ |\x|=1,\quad u(\x,t)=1,\ |\x-\x_{1}|=\epsilon.
\end{equation}
We also take $\Phi_1=1$ and $\ell_1=1$.
The initial concentration is localized at the origin of the sphere
\begin{equation}
u(\x,0)=\Gamma_0\delta(\x),\quad  \Gamma_0 <\frac{4 \pi}{3}.
\end{equation}
The 3D Neumann Green's function in the unit sphere is known explicitly \cite{Cheviakov11}:
\begin{align}
G_0(\x,\hxi)&=\frac{1}{4\pi |\x-\hxi|}+\frac{1}{4\pi |\x|r'} +\frac{1}{4\pi }\ln \left (\frac{2}{1-|\x||\hxi|\cos \theta+|\x|r'}\right ) 
\nonumber \\
&\quad +\frac{1}{6|\Omega|}(|\x|^2+|\hxi|^2)-\frac{7}{10\pi},
\end{align}
where $|\Omega|=4\pi/3$, and
\[\cos \theta =\frac{\x\cdot \hxi}{|\x||\hxi|},\quad \x'=\frac{\x}{|\x|^2},\quad r'=|\x'-\hxi|.
\]
The final constant is chosen so that $\int_{\Omega}G(\x,\hxi)d\x=0$. It follows from equation (\ref{one}) that the accumulation time is
 \begin{align}
T(\x)&= \frac{4\pi/3-\Gamma_0}{4\pi \epsilon D}- \Gamma_0 \left [ G_0(\x,\x_0)-G_0(\x_1,\x_0)\right ]\nonumber \\
&\quad -[4\pi/3-\Gamma_0][G_0(\x,\x_1)-R_0(\x_1,\x_1)]+O(\epsilon).
\end{align}

\begin{figure}[t!]
\centering
\includegraphics[width=13cm]{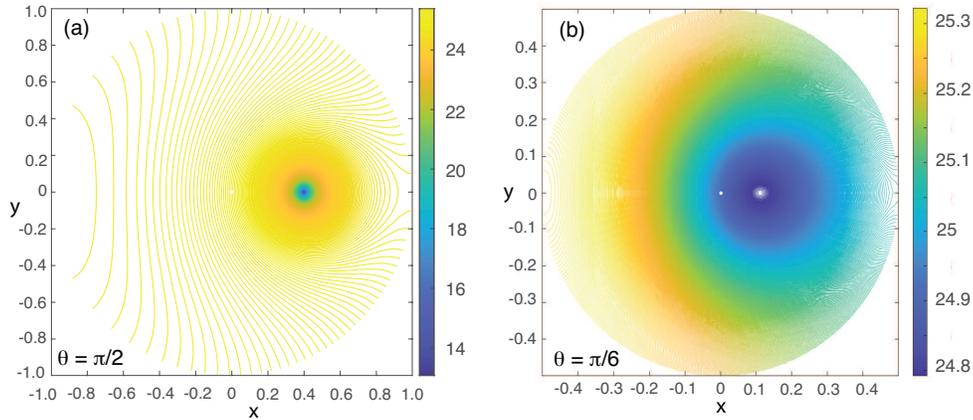} 
\caption{Accumulation time in the unit sphere with a single-hole as shown in Fig. \ref{fig2}(a). Contour plots of the truncated accumulation time $T_0(\x)$ in the $x-y$ plane for $\x=r(\sin \theta \cos \phi,\sin \theta\sin \phi.\cos \theta) $ with $0\leq r \leq 1$, $0\leq \phi <2\pi$ and fixed $\theta$: (a) $\theta=\pi/2$; (b) $\theta=\pi/6$. Other parameter values are the same as in Fig. \ref{fig3}.}
\label{fig4}
\end{figure}

In Fig. \ref{fig3} we show contour plots of the full accumulation time $T(\x)$, $\x=(x,y,z)$, in the $x-y$ plane for several horizontal sections of the sphere (fixed $z$). In the plane containing the initial position $\x_0$ and the center $\x_1$ of the spherical hole, we see that there are two minima of $T(\x)$ located around the points $\x_0$ and $\x_1$, respectively. Note that $T(\x)$ is singular at these points. The singularity in $T(\x)$ as $\x\rightarrow \x_0$ is a consequence of the initial condition involving a Dirac delta function. It is easily removed by taking the initial concentration to be a strongly localized Gaussian, for example. The singularity as $\x\rightarrow \x_1$ is due to the fact that we define $T(\x)$ in terms of the outer solution; it would be resolved by considering the corresponding inner solution. In Fig. \ref{fig4} we present corresponding contour plots for the truncated accumulation time $T_0(\x)$. The reduction in the dependence on the initial position $\x_0$ is clearly seen.

\subsection{Triplet of targets in the unit sphere}

\begin{figure}[t!]
\centering
\includegraphics[width=13cm]{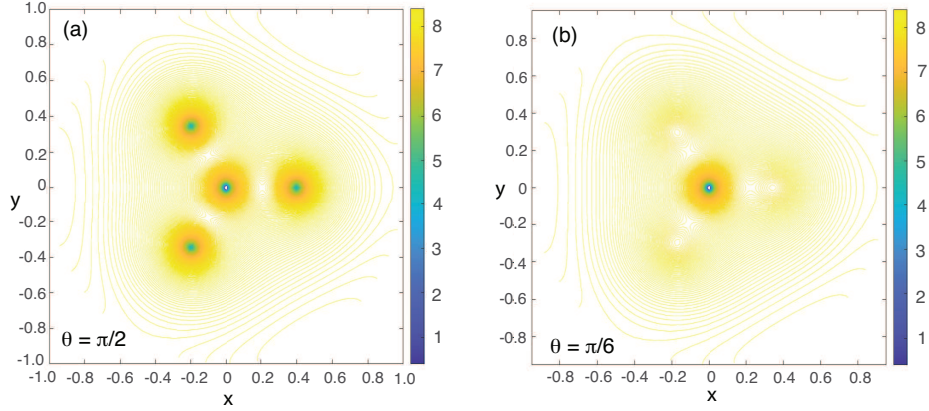} 
\caption{Accumulation time in the unit sphere with a single-hole as shown in Fig. \ref{fig2}(a). Contour plots of the truncated accumulation time $T_0(\x)$ in the $x-y$ plane for $\x=r(\sin \theta \cos \phi,\sin \theta\sin \phi.\cos \theta) $ with $0\leq r \leq 1$, $0\leq \phi <2\pi$ and fixed $\theta$: (a) $\theta=\pi/2$; (b) $\theta=\pi/6$. Other parameter values are the same as in Fig. \ref{fig3}.}
\label{fig5}
\end{figure}

As our second example, consider three identical holes distribution at the points $\x_1=(a,0,0)$, $\x_2 =a(\cos 2\pi/3,\sin 2\pi/3,0)$ and $\x_2 =a(\cos 4\pi/3,\sin 4\pi/3,0)$ in the unit sphere with $a=0.4$. see Fig. \ref{fig2}(b). The boundary conditions are
\begin{equation}
\nabla u(\x,t)\cdot {\bf n}=0,\ |\x|=1,\quad u(\x,t)=1,\ |\x-\x_{1,2,3}|=\epsilon.
\end{equation}
The initial concentration is localized at a point on the $y$-axis so that
\begin{equation}
u(\x,0)=\Gamma_0\delta(\x-\x_0),\quad \x_0=(0,b), \ 0 < b < 1, \Gamma_0 < \pi.
\end{equation}
It follows from equation (\ref{one}) that the accumulation time for three identical targets is
 \begin{align}
T(\x)&= \frac{4\pi/3-\Gamma_0}{12\pi \epsilon D }- \Gamma_0 \left [ G_0(\x,\x_0)-\frac{1}{3}\sum_{j=1}^3G_0(\x_j,\x_0)\right ]
 -\frac{4\pi/3-\Gamma_0}{3}\bigg \{  \sum_{j=1}^3G_0(\x,\x_j) \nonumber \\
&\qquad \qquad -\frac{1}{3}\bigg [\sum_{j=1}^3R_0(\x_j,\x_j)+\sum_{i,j,j\neq i}
G_0(\x_i,\x_j)\bigg \}+O(\epsilon).\nonumber 
\end{align}
In Fig. \ref{fig5} we show contour plots of the $O(1)$ accumulation time $T(\x)$, $\x=(x,y,z)$, in the $x-y$ plane for $\theta=\p/2$ and $\theta=2\pi/5$. As expected, the plots are symmetric with respect to $\phi$-rotations by multiples of $2\pi/3$. There are local minima of $T(\x)$ in the vicinity of the holes and $\x_0$.

\section{Discussion}
In this paper, we continued the development of a method for characterizing the relaxation to a non-trivial steady state of a diffusion process, which is based on the notion of an accumulation time. 
The classical approach is to identify the relaxation rate with the principal non-zero eigenvalue of the negative Laplacian \cite{Ward93,Ward93a}. However, this only yields a global measure of the relaxation rate, and loses all information about the initial position. Moreover, it relies on the existence of a sufficiently large spectral gap. The accumulation time, on the other hand, can be obtained by solving the diffusion equation in Laplace space without any recourse to a spectral decomposition. (One can also consider an eigenfunction expansion of the accumulation time itself, but such an approximation still relies on a spectral gap.) Combining our analysis of diffusion in 3D singularly perturbed domains with our previous study of 2D diffusion \cite{Bressloff22} provides a solid foundation for investigating other relaxation processes in singularly perturbed domains. For example, one could consider more general exterior and interior boundary conditions, provided that there existed a unique steady-state solution. For example, modifying the exterior boundary condition would change the definition of the Green's function used in the outer solution, whereas changing the conditions on the hole boundaries would modify the inner solution and the corresponding singularity conditions for the outer solution. Another generalization, as indicated in section 3(c), would be to consider non-spherical hole shapes, provided that the corresponding shape capacitances could be determined \cite{Ward93,Ward93a}. Finally, one could extend the underlying diffusion equation by including advection terms, for example.

Another class of non-trivial steady state arises within the context of diffusion under stochastic
resetting. The simplest example of such a process is a Brownian particle whose position is reset randomly in time at a constant rate $r$ (Poissonian resetting) to its initial position $\x_0$\cite{Evans11a,Evans11b,Evans14}. One major finding is that the probability density converges to a nonequilibrium stationary state (NESS) that maintains nonzero probability currents. In addition, the approach to the stationary state exhibits a dynamical phase transition, which takes the form of a traveling front that separates spatial regions for which the probability density has relaxed to the NESS from those where it has not. Since the trajectories contributing to the transient region are rare events, one can establish the existence of the phase transition by carrying out an asymptotic expansion of the exact solution \cite{Kus14}. It turns out that this transition can also be understood in terms of the spatial variation of the accumulation time for relaxation to the NESS \cite{Bressloff21C}. That is, $T(\x)\sim |\x-\x_0|/\sqrt{4rD}$ for $|\x-\x_0| \gg \sqrt{D/r}$.

\end{document}